\documentstyle[12pt]{article}

\newcommand{\sect}[1]{\setcounter{equation}{0}\section{#1}}

\newtheorem{theorem}{Theorem}[section]
\newtheorem{proposition}[theorem]{Proposition}
\newtheorem{definition}[theorem]{Definition}

\newtheorem{lemma}[theorem]{Lemma}

\def\be{\begin{equation}}
\def\ene{\end{equation}}
\def\bea{\begin{eqnarray}}
\def\eea{\end{eqnarray}}

\def\fidi{\hskip5pt \vrule height4pt width4pt depth0pt \par}

\def\1{\'{\i}}

\def\al{\alpha}

\def\ve{\varepsilon}

\def\la{\lambda}
\def\si{\sigma}

\def\De{\Delta}

\def\NN{{\bf N}}

\def\ZZ{{\bf Z}}

\def\FSF{{\cal F}}

\def\FSH{{\cal H}}

\def\FSU{{\cal U}}

\def\id{{\rm id}}

\def\Aut{{\rm Aut}\,}

\def\zb{\bar z}
\def\nb{\bar n}
\def\iy{\infty}

\def\tens{\otimes}

\def\indl{{\rm ind}_K^G(\rho_L)}
\def\indr{{\rm ind}_K^G(\rho_R)}

\begin{document}

\title{Induction of quantum group representations}

\author{Nicola Ciccoli \\
Universit\'a degli studi di Perugia \\
Dipartimento di Matematica\\
Via Vanvitelli 1, 06123 Perugia.}

\maketitle
\begin{abstract}
{In this paper we will define a generalized procedure of 
induction of quantum group representations both from quantum
and from coisotropic subgoups proving also their main
properties. We will then show that such a procedure
realizes quantum group representations on generalized quantum
bundles.}
\end{abstract}

\thispagestyle{empty}
{\small {\bf Math. Subj. Classification} 16W30, 17B37, 81R50}


\medskip
{\small {\bf Submitted to:}{\it Journal of Geometry and Physics}}

\sect{Introduction}
\bigskip
Among the various methods to construct representations of a given 
Lie group, the induction procedure plays with no doubt a central role.
On one hand it gives a strong insight in the representation theory
of a large class of non semisimple Lie groups. For example it solves
completely the semidirect product case, linked as it is to the 
orbit method for nilpotent and solvable groups. On the other hand it
gives foundations to the geometric theory of representations of semi
simple Lie groups, from Borel-Weil-Bott theorem on (see for example
\cite{Wal}). Induction of quantum group representations was developed
in \cite{PaWa} and, later, in \cite{GI}. Recently a very detailed work
\cite{GoZa} summarizes results for the compact case. One of the problems
with induced quantum group representations is connected with the extreme
rarity of quantum subgroups. For example in \cite{Gal}, \cite{Mas}, 
\cite{Prz} induced representations for some non semisimple quantum
groups are constructed.  The last two paper, however, does not fit
well the Parshall-Wang setting because they do not start with a well
defined quantum subgroup. To overcome this problem the theory
of coisotropic quantum subgroups introduced in \cite{Ci2} seems
useful.
In this paper we generalize induced quantum group representations
to coisotropic quantum subgroups introduced in \cite{Ci2} and
show how the corresponding corepresentation space of the function
algebra can be intrepreted as the space of sections of an
associated vector bundle to a principal coalgebra quantum bundle,
this last concept being recently introduced by Brzezi\`nski
in \cite{Brz2} and \cite{BrMa2}.
We will not deal with the measure-theoretic concepts needed
to develop correctly all unitariness aspects, which
are the subject of a forthcoming paper.
\bigskip
\bigskip
\sect{Quantum coisotropic subgroups}
\bigskip
In this section we will briefly summarize some basic of the theory
of subgroups of quantum groups, as developed in \cite{Brz1} and
\cite{Ci2} (to which we refer for motivations). A more detailed 
account of the problems we'll deal with can also be found in
\cite{Ci3}, while for the general theory of quantum embeddable 
homogeneous spaces we refer to \cite{DiKo}.

\begin{definition}
Given a Hopf algebra $A_q$ we will call quantum subgroup of $A_q$ any
pair $(B_q,\pi)$ such that $B_q$ is a Hopf algebra and $\pi:A_q\to B_q$
is a Hopf algebra epimorphism. We will call quantum coisotropic left
(resp. right) subgroup any pair $(C_q,\si)$ such that $C_q$ is a coalgebra
and a left (resp. right) $A_q$-module and $\si:A_q\to C_q$ is a
surjective linear map which is also a coalgebra and a left (resp. right)
$A_q$-module homomorphism.
\end{definition}
As usual in the classical case subgroups can be identified through the
kernel of the restriction epimorphism. Thus quantum subgroups are in
1:1 correspondence with Hopf ideals in $A_q$ while left (right)
coisotropic subgroups are in 1:1 correspondece with bilateral coideals
which are also left (right) ideals.

The weaker hypothesis defining coisotropic quantum subgroup is
not descending from a need of generality for its own sake. It is
rather quite naturally imposed by two strictly connected aspects.
First of all, quantum subgroups in the semiclassical limit correspond
to Poisson-Lie subgroups, a fairly rare object. Secondly, while not any
quantum embeddable homogeneous space comes out of a quotient by a
quantum subgroup, a 'big' family of those can be obtained as a quotient
by a coisotropic subgroup. All known families
of embeddable homogeneous spaces are of this kind. Furthermore 
it has been
proven by Etingof and Kazhdan that any quotient by a classical
coisotropic subgroup of a Poisson-Lie subgroup can be quantized.
The same statement does not hold for more general Poisson homogeneous
spaces.

More precisely we have:
\begin{definition}
Let $A_q$ be a Hopf algebra. A right (resp. left) quantum embeddable homogeneous
space of $A_q$ is a subalgebra which is also a right (resp. left) coideal.
\end{definition}
\begin{proposition}
Let $(C,\si)$ be a right (resp. left) coisotropic quantum subgroup of $A_q$.
Then
$$
B_C=\{ f\in A_q \big| (\si\tens id)\De f=\si(1)\tens f\}
$$
is a quantum right embeddable homogeneous space and, respectively
$$
B^C=\{ f\in A_q \big| (id\tens\si)\De f=f\tens\si(1) \}
$$
is a quantum left embeddable homogeneous space.
Conversely let $B_q$ be a right (resp. left) embeddable quantum homogeneous space; then
the right (resp. left) ideal generated by $\{ b-\ve(b)1\}$ is a bilateral coideal
in $A_q$ and identifies a right (resp. left) coisotropic subgroup.
\end{proposition}
The annoying necessity of dealing with both left and right coisotropic
subgroups (which becomes rather involved adding considerations both
on left and right corepresentations) disappear when considering the
unitary side of the theory, which forces precise choices.

\bigskip
\bigskip
\sect{Quantum induction}
\bigskip
In this section we define induced corepresentations from coisotropic
quantum subgroups and prove their main properties.
Let us begin recalling some definitions. Let $A$ be a Hopf algebra and let 
$V$ be a vector space (we're considering a fixed base field of zero 
charachteristic). $V$ is said a right corepresentation for $A$ if there
exists a linear map
\be
\psi_R:V\rightarrow V\tens A
\ene
such that:
\be
\label{coinv}
(id\tens \ve)\circ\psi_R=id_V
\ene
\be
\label{cove}
(id\tens \Delta)\circ\psi_R=(\psi_R\tens\id_A )\circ\psi_R
\ene
We will say that $V$ is a left corepresentation for $A$ if there exists
a linear map 

\be
\psi_L:V\rightarrow A\tens V
\ene
such that:
\be
(\ve\tens id)\circ\psi_L=id_V
\ene
\be
(\Delta\tens id)\circ\psi_L=(id_A\tens\psi_L )\circ\psi_L
\ene
Let us remark explicitly that above definitions are well posed for
any coalgebra $A$. 

Let us consider now fixed a quantum group $A_q=A_q(G)$ and a coisotropic
quantum subgroup $(C_q,\pi)=(A_q(K),\pi)$ (when we do not specify coisotropic
subgroups to be left or right we mean the result is valid in both cases).

\begin{proposition}
\label{rappre}
The map $R=(id\tens\pi)\circ\De_G:A_q(G)\to A_q(G)\tens A_q(K)$ defines a
right corepresentation of $A_q(K)$ on $A_q(G)$. Similarly the map
$L=(\pi\tens id)\circ\De_G$ defines a left corepresentation.
\end{proposition}
Proof:

We will verify conditions \ref{coinv} and \ref{cove} for $R$. The proof
for $L$ proceeds along the same lines. Let us start with \ref{cove}:
$$
Rf=\sum\limits_{(f)} f_{(1)}\tens\pi(f_{(2)})
$$
$$
(id\tens \De_K)(Rf)=\sum\limits_{(f)} f_{(1)}\tens(\De_K\circ\pi)(f_{(2)})=
\sum\limits_{(f)}f_{(1)}\tens \pi(f_{(2)})\tens \pi(f_{(3)}) 
$$
where we have used coassociativity of $\De_G$ together with the fact that
$\pi$ is a coalgebra morphism. On the other side
$$
((R\tens id)\circ R)f=\sum\limits_{(f)} R(f_{(1)})\tens \pi (f_{(2)})=
\sum\limits_{(f)}(id\tens \pi)(f_{(1)_{(1)}}\tens f_{(1)_{(2)}})\tens
\pi (f_{(2)}) =$$
$$
= \sum\limits_{(f)}f_{(1)}\tens\pi(f_{(2)})\tens\pi(f_{(3)})
$$
again using coassociativity; this proves \ref{cove}. For what concerns
\ref{coinv} we have:
$$
(id\tens\ve_K)(Rf)=\sum\limits_{(f)}f_{(1)}\tens(\ve_K\circ\pi)(f_{(2)})=
$$
$$
=\sum\limits_{(f)}f_{(1)}\tens\ve_G(f_{(2)})=((id\tens\ve_G)\circ\De)f=f
$$
The proof for $L$ goes along in the same way.
\fidi
\begin{lemma}
The following identities hold:
\be
(\De_G\tens id )\circ R=(id\tens R)\circ\De_G
\ene
\be
(id\tens \De_G)\circ L=(L\tens id)\circ\De_G
\ene
\end{lemma}
Proof

Let us prove the first, the second is proven in the same way.
$$
(\De_G\tens id)(Rf)=\sum\limits_{(f)}\De_G(f_{(1)})\tens\pi(f_{(2)})=
\sum\limits_{(f)}f_{(1)}\tens f_{(2)}\tens\pi(f_{(3)})=
$$
$$
(id\tens R)(\De_G f)=\sum\limits_{(f)}f_{(1)}\tens R(f_{(2)})=\sum\limits
_{(f)} f_{(1)}\tens f_{(2)}\tens \pi(f_{(3)})
$$
\fidi
Straightforward calculations allow to prove also the following claim:
\begin{lemma}
Let $(A_q(K),\pi)$ be a left coisotropic subgroup; then we the 
multiplication properties:
\begin{equation}
L(fg)=\De_G(f)L(g); \qquad R(fg)=\De f \cdot R(g) 
\qquad \forall f, g\in\FSF_q(G)
\end{equation}
where $\cdot$ denotes the action of $\FSF_q(G)\tens\FSF_q(G)$ on
$\FSF_q(G)\tens A_q(K)$ given by $(f\tens g)\cdot(h\tens c)=
fh\tens g\cdot c$. Similarly if $(A_q(K),\pi)$ is a right coisotropic
subgroup we have:
\begin{equation}
L(fg)=L(f)\De_G(g);\qquad R(fg)=R(f)\De_G(g) \qquad \forall f,g\in\FSF_q(G)
\end{equation}
\end{lemma}

Let us start now from a right corepresentation $\rho_R$ and a left
corepresentation $\rho_L$ of the coalgebra $A_q(K)$ on the vector 
space $V$. Define:
\be
\label{leftind}
{\rm ind}_K^G(\rho_L)=\{ F\in A_q(G)\tens V\big| \quad (R\tens id)(F)=
(id\tens\rho_L)(F)\}
\ene
\be
\label{rightind}
{\rm ind}_K^G(\rho_R)=\{ F\in V\tens A_q(G)\big| \quad (id\tens L)F=
(\rho_R \tens id)(F)\}
\ene
Both these spaces are kernels of certain linear operators and thus are
closed vector subspaces of $A_q(G)\tens V$ and $V\tens A_q(G)$.
\begin{proposition}
$\De_G\tens id$ defines a left corepresentation of $A_q(G)$ on ${\rm ind}_K^G
(\rho_L)$ and $id\tens \De_G$ defines a right corepresentation of the 
same space on ${\rm ind}_K^G(\rho_R)$.
\end{proposition}
Proof

Let us prove the first claim. To begin with we need to prove that
$(\De_G\tens id)(F)\in A_q(G)\tens\indl$. Linearity of all the operations
involved implies that we can limit ourselves to the case $F=f\tens v$.
Then what we have to prove is equivalent to
$$
(id\tens R\tens id)(\De_G\tens id)(f\tens v)=(id\tens id\tens\rho_L)
(\De_G\tens id)(f\tens v)
$$
Using formula 3.7 in lemma 3.1 the first term is equal to
$$
((id\tens R)\circ\De_G)\tens id)(f\tens v)=((\De_G\tens id)\circ R)
\tens id)(f\tens v)=(\De_G\tens id\tens id)(R\tens id)(f\tens v)
$$
which, applying the fact that $f\tens v$ is in $\indl$, is equal to
$$
=(\De_G\tens id\tens id)(id\tens\rho_L)(f\tens v)=(id\tens id\tens\rho_L)
(\De_G\tens id)(f\tens v)
$$
Now the corepresentation condition \ref{cove} is equivalent to coassociativity
of $\De_G$ and the condition \ref{coinv} is equivalent to properties
of the counity $\ve_G$. The second claim follows similarly from 3.8.
\fidi
\begin{definition}
Given a right corepresentation $\rho_R$ of the quantum coisotropic 
subgroup $(A_q(K),\pi)$
the corresponding corepresentation $\indr$ of $A_q(G)$ is called
induced representation.
\end{definition}
{\bf Remark}
In case $\rho_L$ is a one-dimensional corepresentation the induced 
representation is given by
$$
\indl =\{ f\in A_q(G)\big| (R\tens id)f=(id\tens\rho_L)f\}
$$
with coaction $\De_G$. Similarly if $\rho_R$ is one-dimensional 
$$
\indr =\{ f\in A_q(G)\big| (id\tens L)f=(\rho_R\tens id)f\}
$$
In both cases they can be seen as subrepresentation of the regular
representation (left or right) $(A_q(G),\De_G)$. In analogy with classical
terminology we will call such representations {\it monomial}.

As usual in representation theory we will often idntify the space of the
representation with the representation itself, being clear which is 
the representation map.

\begin{proposition}
If $\rho_R$ and $\rho'_R$ are equivalent right corepresentations of 
$A_q(K)$ then $\indr$ and ${\rm ind}_K^G(\rho'_R)$ are equivalent right
corepresentations. An obvious statement is valid for the case of left
corepresentations.
\end{proposition}
Proof

Let $\rho_R:V \to A_q(K)\tens V$ and $\rho'_R:W\to A_q(K)\tens W$ be equivalent.
Then there exists a vector space isomorphism $F:V\to W$ such that
$\rho'_R\circ F=(id\tens F)\circ\rho_R$. Define then the isomorphism
$$
{\tilde F}=F\tens id:V\tens A_q(G)\to W\tens A_q(G)
$$
Let us prove that $ {\tilde F}(\indr)\subset {\rm ind}_K^G(\rho'_R)$. Let
$v\tens f\in \indr$. We want to prove that $F(v)\tens f$ verifies
$$
F(v)\tens L(f)=\rho'_R(F(v))\tens f
$$
We have that
$$
\rho'_R(F(v))\tens f=(id\tens F)\circ\rho_R)(v)\tens f=
(id\tens F\tens id)(\rho_R(v)\tens f)=
$$
$$
=(F\tens id\tens id)(v\tens L(f))=F(v)\tens L(f)
$$
Let us prove now that the vector space isomorphism ${\tilde F}$ intertwines
corepresentations, i.e.
$$
(id\tens \De_G)\circ{\tilde F}=({\tilde F}\tens id)(id\tens \De_G)
$$
We have:
$$(id\tens \De_G)\circ{\tilde F}(v\tens f)=(id\tens \De_G)(F(v)\tens f)=
F(v)\tens\De_Gf=
$$
$$
=({\tilde F}\tens id)(v\tens \De_Gf)=({\tilde F}\tens id)(id\tens \De_G)
(v\tens f)
$$
\fidi
Another major property is the behaviour of the induction procedure with 
respect to direct sum. Let us recall that if $\rho_1$ and $\rho_2$ are
right corepresentations of the coalgebra $C$ on the vector spaces $V$
and $W$ respectively, then we define:
\be
\label{sum}
\rho_1\oplus\rho_2:V\bigoplus W\to C\tens (V\bigoplus W)
\ene
$$
\rho_1\oplus \rho_2=(id\tens \imath_V)\circ \rho_1\circ p_V+(id\tens
\imath_W)\circ\rho_2\circ p_W
$$
where $\imath_V:V\to V\bigoplus W$ and $\imath_W:W\to V\bigoplus W$ are
the natural immersion and $p_V:V\bigoplus W\to V$ and $p_W:V\bigoplus W
\to W$ are the natural projections. 
\begin{proposition}
\label{sumind}
Let $(\rho_k,V_k)$, $k\in\NN$ be right $A_q(K)$-corepresentations of a given
coisotropic quantum subgroup $(A_q(K),\pi)$ of $A_q(G)$. Then we have
an equivalence of right corepresentations
$$
{\rm ind}_K^G(\bigoplus_k \rho_k)\equiv
\bigoplus_k {\rm ind}_K^G(\rho_k)
$$
In particular if ${\rm ind}_K^G(\rho)$ is irreducible then $\rho$ is
irreducible.
\end{proposition}
Proof

We will prove the proposition for a finite sum.
Due to associativity of the direct sum it is sufficient to prove the
theorem for $k=1,2$. We will begin proving that $p_V\tens id+p_W\tens id$
is an isomorphism of vector spaces of ${\rm ind}_K^G(\rho_1\oplus \rho_2)$
in ${\rm ind}_K^G(\rho_1)\bigoplus {\rm ind}_K^G(\rho_2)$.
Due to the fact that the above map is a linear isomorphism of 
$(V\bigoplus W)\tens A_q(G)$ in $(V\tens A_q(G))\bigoplus (W\tens A_q(G))$
it is sufficient to prove that $(p_V\tens id)F\in {\rm ind}_K^G(\rho_1)$
and $(p_W\tens id)F\in {\rm ind}_K^G(\rho_2)$.
We have:
$$
(id\tens L)(p_V\tens id)F=(p_V\tens id\tens id)(id\tens L)F=
(p_v\tens id\tens id)((\rho_1\oplus \rho_2)\tens id)F=
$$
$$
=((p_V\tens id)\circ (\rho_1\oplus \rho_2))\tens id (F)=
((p_V\tens id)[(\imath_V\tens id)\circ\rho_1\circ p_V+(\imath_W\tens id)
\circ\rho_2\circ p_W]\tens id(F)=
$$
$$
=[(\imath_V\tens id)\circ \rho_1\circ p_V\tens id](F)=(\rho_1\tens id)(
p_V\tens id)F
$$
and similarly for $(p_W\tens id)$. We can then identify the vector spaces
writing for every $f\in {\rm ind}_K^G(\rho_1\oplus\rho_2)$, $f=f_1+f_2$ with
$f_1\in{\rm ind}_K^G(\rho_1)$ and $f_2\in{\rm ind}_K^G(\rho_2)$.
Let us prove then the intertwining property:
$$
(id\tens \De_G)f=(id\tens \De_G)(f_1+f_2)=(id\tens \De_G)f_1+(id\tens \De_g)
f_2$$
which is the sum of the induced corepresentations.
\fidi
Another interesting properties of the classical induction procedure
that we want to mimic is the so called double induction.
At this purpose let us remark that a coisotropic quantum subgroup cannot
have a quantum subgroup but only coisotropic quantum subgroups. Let us
consider then the case in which $A_q(K)$ is a quantum subgroup of $A_q(G)$,
coisotropic or not, and $A_q(H)$ is a coisotropic quantum subgroup of
$A_q(K)$.

\begin{proposition}
Let $\rho_R$ be a right corepresentation of $A_q(H)$ on $V$. Then there is an
equivalence of right $A_q(G)$-corepresentations between
$$
{\rm ind}_H^G(\rho_R)\equiv {\rm ind}_K^G({\rm ind}_H^K(\rho_R))  
$$
The same property holds for left corepresentations.
\end{proposition}

Proof

Let us denote with $\pi_{KH}$ and $\pi_{GK}$ respectively the maps defining
$H$ as a quantum subgroup of $K$ and $K$ as a quantum subgroupof $G$ and let
$\pi_{GH}=\pi_{KH}\circ\pi_{GH}$ defining $H$as a quantum subgroup of $G$.
Let us denote with $L_{GK}$, $L_{KH}$ and $L_{GH}$ the corresponding
left corepresentations granted by \ref{rappre}. Let us remark that $L_{GH}=
(\pi_{KH}\tens id)\circ L_{GK}$. The required isomorphism between
representation spaces is given by:
$$
id\tens L_{GK}:{\rm ind}_H^G(\rho)\to {\rm ind}_K^G({\rm ind}_H^K(\rho))  
$$

To prove it remark that elements of the second space are those $v\tens g
\tens f$ in $V\tens A_q(K)\tens A_q(G)$ verifying
$$
1)\qquad \rho_R(v)\tens g\tens f=\sum_{(g)} v\tens\pi_{KH}(g_{(1)})
\tens g_{(2)}\tens f
$$
$$
2)\qquad v\tens\De_K(g)\tens f=\sum_{(f)}v\tens g\tens\pi_{GK}(f_{(1)})
\tens f_{(2)}
$$
Applying $id\tens\ve_K\tens id\tens id$ to this second equality gives
$$
v\tens g\tens f=\sum_{(f)}\ve_K(g) v\tens\pi_{GK}(f_{(1)})\tens f_{(2)}
$$
from which one can prove that $v\tens g\tens f\mapsto \ve(g)v\tens f$
is both a left and a right inverse of $id\tens L_{GK}$. The fact
that the actions are intertwined results from coassociativity of $\De_G$.
\fidi

Let us prove now how the induction procedure interacts with automorphisms
of the Hopf algebra structure.
Let $A_q(G)$ be the quantum group and $(A_q(K),\pi)$ a coisotropic quantum
subgroup. Let $\al\in \Aut (A_q(G))$ be a Hopf-algebra automorphism. Then
there exists one and only one coalgebra automorphism ${\tilde \al}
:A_q(K)\to A_q(K)$ such that: ${\tilde \al}\circ\pi=\pi\circ\al$.

\begin{proposition}
In the above hypothesis we have the equivalence of representations
$$
{\rm ind}_K^G((id\tens{\tilde \al})\circ\rho_R)=(id\tens\al)({\rm ind}_K^G
(\rho_R))
$$
\end{proposition}

Proof

The vector space of the corepresentation on the left is:
$$
\{ F\in V\tens A_q(G)\big|\qquad (id\tens L)F=(id\tens{\tilde \al}\tens id)
(\rho_R\tens id)F \}
$$
If $F$ belongs to ${\rm ind}_K^G(\rho)$ then
$$
(id\tens L)(id\tens\al)F=v\tens L(\al(f))=v\tens((\pi\tens id)\De(\al (f))=
$$
$$
=v\tens((\pi\tens id)(\al\tens\al(\De f)))=v\tens(({\tilde \al}\tens{\tilde
\al})(\pi\tens id(\De f)))=(id\tens{\tilde \al}\tens{\tilde \al})(v\tens L(f))
$$
$$
=(id\tens{\tilde \al}\tens{\tilde \al})(v\tens\rho(f))=(id\tens{\tilde \al})
(\rho(v))\tens \al(f)
$$
and this proves the isomorphism between the vector spaces carrying the
representation. The intertwining property is a simple consequence of 
$\al$ being a coalgebra morphism.
\fidi

\bigskip
\bigskip
\sect{Geometric realization on homogeneous quantum bundles}
\bigskip
The purpose of this section is to explicit the
relations between induced corepresentations from coisotropic subgroups
and embeddable quantum homogeneous spaces. We will split the
two cases of quantum and coisotropic subgroup, although the first
one is a special case of the second, to point out the differences.
Much of what follows is both inspired from and intimately related to 
results in references \cite{Brz2}, \cite{BrMa1}, \cite{BrMa2}, where
the theory of quantum principal bundles with structure group given
by a coalgebra has been developed. Those results are reinterpreted
(and slightly generalized) here in the context of induced corepresentations.
This reflects what happens in the classical case where there is
a bijective correspondence between induced representations and
homogeneous vector bundles (i.e. vector bundles associated to
principal bundles, \cite{Wal}). For a more detailed relation
on the quantum bundle interpretation we refer to \cite{GoZa}, which
appeared while this paper was in preparation and which contains
more details, although limited to the compact and subgroup case.

Let $\FSF_q(G)$ be a Hopf algebra with invertible antipode and
let $(\FSF_q(K),\pi)$ be a quantum subgroup. Let $B_K$ be the
corresponding quantum quotient space
$$
B_K=\{ f\in\FSF_q(G) \quad \big| \quad(\pi\tens id)\De(f)=1\tens f\}
$$
Let $\rho_R$ be a $\FSF_q(K)$ right corepresentation on $V$
and let $\FSH_R$ be the space of the induced $\FSF_q(G)$ corepresentation.

\begin{lemma}
$\FSH_R$ is a left and right $B_K$-module, where the action is given
by the linear extension of
$$
b\cdot(f\tens v)=(bf\tens v)
$$
$$
(f\tens v)\cdot b=fb\tens v
$$
\end{lemma}

Proof

Let $v\tens f\in\FSH_R$. Then
$$
L(v\tens bf)=\sum_{(b)(f)}v\tens \pi(b_{(1)}f_{(1)})\tens b_{(2)}f_{(2)}
=\sum_{(b)(f)}v\tens \ve(b_{(1)})\pi(f_{(1)})\tens b_{(2)}\pi(f_{(2)})=
$$
$$
=\sum_{(f)}v\tens f_{(1)}\tens bf_{(2)}=\sum_{(v)}v_{(0)}\tens v_{(1)}
\tens bf=b\cdot (\rho_R\tens id)(v\tens f)
$$
The same proof holds for right action ($\pi$ is an algebra morphism).
\fidi
Let us now recall that linear maps $f,g:\FSF_q(K)\to\FSF_q(G)$
can be multiplied according to convolution:
$$
(f\ast g)(k)=\sum_{(k)}f(k_{(1)})g(k_{(2)})
$$
If $f:\FSF_q(K)\to\FSF_q(G)$ its convolution inverse, if it exists,
is a map $f^{-1}:\FSF_q(K)\to\FSF_q(G)$ such that:
$$
\sum_{(k)}f(k_{(1)})f^{-1}(k_{(2)})=\ve(k)1=\sum_{(k)}f^{-1}
(k_{(1)})f(k_{(2)})
$$
\begin{definition}
\label{section1}
A quantum subgroup is said to have a left (resp. right) section if there exists a linear,
convolution invertible, map $\phi:\FSF_q(K)\to\FSF_q(G)$ such that:

i) $\phi(1)=1$;

ii) $(\pi\tens id)\De_G\circ\phi=(1\tens \phi)\circ\De_K$ (or, respectively
ii') $(id\tens\pi)\De_G\circ\phi=(\phi\tens 1)\circ\De_K$.

If, furthermore, $\phi$ is an algebra morphism then $(\FSF_q(K),\pi)$
is said to be trivializable.
\end{definition}
The second condition is an intertwining condition between the corepresentation
of $\FSF_q(K)$ on itself and its corepresentation on $\FSF_q(G)$.
\begin{proposition}
\label{iso-hom1}
Let us suppose that $(\FSF_q(K),\pi)$ has a section $\phi$. Then:

i) $\FSF_q(G)$ is isomorphic to $\FSF_q(K)\tens B_K$ as a vector space;

ii) $\FSH_R$ is isomorphic to $V\tens B_K$ as a $B_K$-module (both left and
right).
\end{proposition}

Proof

Let us consider, first of all, the following lemma, which can be proved
with usual Hopf algebra techniques.
\begin{lemma}
The convolution inverse of the section $\phi$ verifies
$$
(\pi\tens id)\De_G\phi^{-1}=(S\tens \phi^{-1})\tau_{1,2}\De_K
$$
where $\tau_{1,2}$ is the map interchanging terms in tensor product.
\end{lemma}

Now let us consider $f\in\FSF_q(G)$. Then
$$
f=\sum_{(f)}\ve(f_{(1)})f_{(2)}=\sum_{(f)}\ve(\pi(f_{(1)}))f_{(2)}
=\sum_{(f)}\phi(\pi(f_{(1)}))\phi^{-1}(\pi(f_{(2)}))f_{(3)}
$$
Let us now prove that $\sum_{(f)}\phi^{-1}(\pi(f_{(1)}))f_{(2)}$
belongs to $B_K$.

$$
(\pi\tens id)\De(\sum_{(f)}\phi^{-1}(\pi(f_{(1)}))f_{(2)})=
\sum_{(f)}\pi(\phi^{-1}(\pi(f_{(1)}))_{(1)})f_{(2)_{(1)}}\tens
\phi^{-1}(\pi(f_{(1)}))_{(2)}f_{(2)}))_{(2)}=
$$
$$
=\sum_{(f)}S_K(\pi(f_{(2)}))\pi(f_{(3)})\tens\phi^{-1}(\pi(f_{(1)}))f_{(4)}
$$
$$
=\sum_{(f)}1\tens \ve(f_{(2)})\phi^{-1}(\pi(f_{(1)}))f_{(3)}=
$$
$$
=1\tens\sum_{(f)}\phi^{-1}(\pi(f_{(1)}))f_{(2)}
$$
This proves that the map
$A_{\phi}:\FSF_q(K)\tens B_K\to\FSF_q(G)$ linearly extending 
$k\tens b\mapsto \phi(k) b$
is surjective and its right inverse is given by 
$$
A_{\phi}^{-1}:f\mapsto \sum_{(f)} \pi(f_{(1)})\tens\phi^{-1}(\pi(f_{(2)}))
f_{(3)}
$$
Let us prove that $A_{\phi}^{-1}$ is also a left inverse:
$$
A_{\phi}^{-1}(\phi(k)b)=\sum_{(b)(\phi(k))}\pi(\phi(k)_{(1)}b_{(1)})\tens
\phi^{-1}(\pi(\phi(k)_{(2)}b_{(2)}))\phi(k)_{(3)}b_{(3)}=
$$
using the fact that $ii)$ of definition 4.1 implies that $(\pi\tens\pi\tens id)
(\De_G\tens id)\De_G\phi$ equals $(id\tens id\tens\phi)(\De_K\tens id)
\De_K$ we have
$$
=\sum_{(\phi(k))}\pi(\phi(k)_{(1)})\tens\phi^{-1}(\pi(\phi(k)_{(2)}))
\phi(k)_{(3)}b=\sum_{(k)}k_{(1)}\tens\phi^{-1}(k_{(2)})
\phi( k_{(3)})b
$$
$$
=\sum_{(k)}\ve(k_{(2)})k_{(1)}\tens b=k\tens b
$$
which proves the claim.

The proof of the second isomorphism, although more involved, is
quite similar to the preceding one. Let us define first of all
$$
T_{\phi}:V\to\FSH_R;\qquad v\mapsto\sum_{(v)}v_{(0)}\tens\phi(v_{(1)})
$$
We have to prove that $T_{\phi}$ really takes its values in $\FSH_R$.
This is true if and only if
$$
\sum_{(v)}\rho_R(v_{(0)})\tens\phi(v_{(1)})=\sum_{(v)}v_{(0)}\tens
(\pi\tens id)\De_G(\phi(v_{(1)}))
$$
which is a straightforward calculation. Let us define now
$I_{\phi}:V\tens B_K\to\FSH_R$ as the linear extension of
$v\tens b\to T_{\phi}(v)\cdot b$.
We want to prove that $I_{\phi}$ is bijective, which we will do by
showing that its bilateral inverse is the linear extension of
$$
I^{\phi}(v\tens f)=\sum_{(f)}v\tens \phi^{-1}(\pi(f_{(1)}))f_{(2)}=
\sum_{(v)}v_{(0)}\tens \phi^{-1}(v_{(1)})f
$$
First of all
$$
I^{\phi}(I_{\phi}(v\tens b))=\sum_{(v)}I^{\phi}(v_{(0)}\tens \phi(v_{(1)}))b
$$

$$
=\sum_{(v)(b)(\phi(v_{(1)}))}v_{(0)}\tens \phi^{-1}(\pi(\phi(v_{(1)})_{(1)}
b_{(1)}))\phi(v_{(1)})_{(2)}b_{(2)}
$$
$$
=\sum_{(v)(\phi(v_{(1)}))}v_{(0)}\tens\phi^{-1}(\pi(\phi(v
_{(1)}))_{(1)})\phi(v_{(1)})_{(2)}b=\sum_{(v)}v_{(0)}\tens\phi^{-1}
(v_{(1)_{(1)}})\phi(v_{(1)_{(2)}})b=v\tens b
$$
The fact that $I^{\phi}$ takes its values in $V\tens B_K$ follows
from the first part of the proof. Finally
$$
I_{\phi}(\sum_{(f)}v\tens \phi^{-1}(\pi(f_{(1)}))f_{(2)})
=I_{\phi}(\sum_{(v)}v_{(0)}\tens \phi^{-1}(v_{(1)})f)=
$$
$$
=\sum_{(v)}v_{(0)}\tens \phi(v_{(1)_{(1)}})\phi^{-1}(v_{(1)_{(2)}})f=
v\tens f
$$
\fidi

The second isomorphism of the proposition can be considered a morphism
of $\FSF_q(G)$-comodules provided we put on $V\tens B_K$ the right 
coaction. Explicitely it is given by:
$$
\si_R:V\tens B_K\to V\tens B_K\tens \FSF_q(G)
$$
$$
\si_R(v\tens b)=\sum_{(b)(\phi(v_{(2)}))(v)}v_{(0)}\phi^{-1}(v_{(0)_{(1)}})
\phi(v_{(1)})_{(1)}\tens b_{(1)}\tens\phi(v_{(1)})_{(2)}b_{(2)}
$$
$$
=\sum_{(v)(b)} v_{(0)}\tens b_{(1)}\tens \phi(v_{(1)})b_{(2)}
$$
Remark also that in the case of {\it right sections}
we would have obtained a $B^K$-module homomophism
for $\FSH^L$ both on left and right side.

Let us now consider the
more delicate situation in which we start from a right corepresentation 
$\rho_R$ of a left (resp. right) coisotropic quantum subgroup $(C,\pi)$ of 
$\FSF_q(G)$. Let $B_C$ be the corresponding embeddable quantum
homogeneous space. Then one can multiply functions on the space of
induced representation by functions on the homogeneous space only on
one side.
\begin{lemma}
If $(\FSF_q(K),\pi)$ is a left (resp. right) coisotropic quantum subgroup then 
$\FSH_R$ is a right (resp. left) sub-$B_C$-module (resp. sub $B^C$-module)
of $V\tens\FSF_q(G)$.
\end{lemma}
Proof

Let us consider the case of a right coisotropic subgroup.
$$v\tens(\pi\tens id)\De(fb)=\sum_{(f)(b)}v\tens\pi(f_{(1)}b_{(1)})\tens
f_{(2)}b_{(2)}
$$
$$
=\sum_{(f)(b)}v\tens \pi(f_{(1)}\ve(b_{(2)})b_{(1)})\tens f_{(2)})=
\sum_{(v)}v_{(0)}\tens v_{(1)}\tens fb
$$
\fidi

Let us remark that the convolution product
of linear maps $f, g:C\to\FSF_q(G)$ is still well defined.
\begin{definition}
A right (resp. left) coisotropic subgroup $(C,\pi)$ is said to have
a section $\phi$ if there exists a linear map $\phi:C\to\FSF_q(G)$
convolution invertible and such that:

i) $\phi(\pi(1))=1$;

ii) $\sum_{(c)}\pi(\phi(c)_{(1)}u)\tens\phi(c)_{(2)}=\sum_{(v)}
\pi(v_{(1)}u)\tens\phi(\pi(v_{(2)}))$ for every $c\in C$,
$u\in\FSF_q(G)$ and $v\in\pi^{-1}(c)$. 

(resp. ii') $\sum_{(c)}\phi(c)_{(1)}\tens\pi(u\phi(c)_{(2)})=\sum_{(v)}
\phi(\pi(v_{(1)})\tens\pi(uv_{(2)})$ for every $c\in C$, $v\in\FSF_q(G)$
and $v\in\pi^{-1}(c)$.

If, furthermore, $\phi$ is a right (resp. left) $\FSF_q(G)$-module map then
it is said to be a trivialization.
\end{definition}
\begin{lemma}
\label{phimenuno}
A right coisotropic subgroup section verifies
$$
\sum_{(c)}\pi(\phi^{-1}(c)_{(1)}u)\tens\phi^{-1}(c)_{(2)}=
\sum_{(v)}\pi(S(v_{(2)}u)\phi^{-1}(\pi(v_{(1)}))
$$
where $c\in C$, $v\in\pi^{-1}(c)$.
\end{lemma}

Let us remark that $u=1$ yields again condition (ii) of \ref{section1}.

\begin{proposition}
\label{iso-hom2}
Let us suppose that $(C,\pi)$ is a right coisotropic subgroup
with a section $\phi$. Then:

i) $\FSF_q(G)$ is isomorphic to $C\tens B_C$ as a vector space;

ii) $\FSH_R$ is isomorphic to $V\tens B_C$ as a left $B_C$-module.
\end{proposition}
Proof

Let $f\in \FSF_q(G)$. Using lemma \ref{phimenuno} one can prove like
in \ref{iso-hom1} that $\sum_{(f)}\phi^{-1}(\pi(f_{(1)}))f_{(2)}$
belongs to $B_C$. The inverse isomorphisms $A_{\phi}$ and $A_{\phi^{-1}}$
are then realized exactly as in \ref{phimenuno}. Also the second part
of the proof goes along exactly in the same way.
\fidi

Remark that the $B_C$-module isomorphism in \ref{iso-hom1} and \ref{iso-hom2}
proves that $\FSH_R$ is a projective $B_C$-module, which is a natural
property to ask, as generalization of Swan's theorem, to spaces of
sections (see \cite{Con}).

The explicit isomorphism of \ref{iso-hom2} can be used to describe the
corepresentation directly on this space. The 
coaction is given by:

$$
\rho(b\tens v)=\sum_{(b)(v)(v_{(1)})(\phi(v_{(2)}))}
b_{(1)}\phi (v_{(2)})_{(1)}\tens b_{(2)}\phi(v_{(2)})_{(2)}\phi^{-1}
(v_{(1)})\tens v_{(0)}
$$
$$
=\sum_{(v)(b)}v_{(0)}\tens b_{(1)}\tens\phi(v_{(1)})b_{(2)}
$$
this allow a complete description of the induced corepresentation
from the following data: the homogeneous space $B_C$ and the section
$\phi$. It is very useful in applications as the following examples
show.

Let us give two examples in which one cannot avoid the use of coisotropic 
subgroups. We start with the non standard Euclidean quantum group $E_{\kappa}(2)$
as described in \cite{Dolmo}. The coisotropic subgroup we will
consider is the coalgebra $C$ with a denumerable family of group-like
generators $c_p, p\in\ZZ$ and with restriction epimorphism
$$
\pi(v^p)=c_p\qquad \pi(a_1)=\pi(a_2)=0
$$

The corresponding quantum homogeneous space is the $\kappa$-plane
$$
[a_1,a_2]=\kappa (a_1-a_2)
$$
As shown in \cite{Brz3} this quantum homogeneous space has a section
$$
\phi:C\to E_{\kappa}(2);\qquad c_p\mapsto v^p
$$
Starting from one-dimensional irreducible corepresentations of this
subgroup $\rho_n:1\to 1\tens c_n$ we obtain as corepresentation space the
$\kappa$-plane with coaction
$$
{\it ind}(\rho_n)(a_1)=(v+v^{-1})v^n\tens a_1-i(v-v^{-1})v^n\tens a_2+a_1v^n\tens 1
$$
$$
{\it ind}(\rho_n)(a_2)=i(v-v^{-1})v^n\tens a_1+(v+v^{-1})v^n\tens a_2+a_2v^n\tens 1
$$
For $n=0$ this is simply the regular corepresentation on the $\kappa$-
plane. Using the duality pairing explicitely described in \cite{Dolmo}
it is easy to derive an algebra representation of the corresponding
quantized enveloping algebra, which is never irreducible. The corresponding
decomposition into irreducibles amounts to $E_{\kappa}(2)$-harmonic
analysis on the $\kappa$ plane and has been carried through in 
\cite{Dolmo} for the $n=0$ case.

The second, similar, example is given by the standard Euclidean quantum group
$E_q(2)$ with the family of right coisotropic subgroups $(C,\pi_{\la})$,
where $C$ is as before and the restriction epimorphism is given by
$$
\pi(v^rn^s{\bar n}^k)=q^{2r(k+s)}\la^s{\bar\la}^kc_r(c_{-1};q^{-2})_k
(c_1;q^2)_s
$$
The corresponding quantum homogeneous spaces are called quantum hyperboloids
and described in \cite{BCGST}and \cite{Ci2}; they are algebras in two generators, $z$ and 
$\zb$, with relation $z\zb=q^2\zb z+(1-q^2)$. 
In \cite{Brz3} these subgoups are shown to admit sections
$$ \phi_r:C\to E_q(2);\qquad c_p\mapsto v^{p-r}$$
Starting from one-dimensional irreducible corepresentations $\rho_m:1\mapsto 1\tens
c_m$ we obtain corepresentations on the quantum hyperboloids given by
$$
{\it Ind}(\rho_m)(z)=v^{m-r+1}\tens z+v^m n\tens 1
$$
$$
{\it Ind}(\rho_m)({\bar z})=v^{m+r-1}\tens\zb+v^m\nb\tens 1
$$
The decomposition into irreducibles for the $\iy$-dimensional corresponding
$\FSU_q(e(2))$-representation is dealt with, in the $r=0$-case, in
\cite{BCGST}.
\bigskip\bigskip
\section{Acknowledgements}
The author was partly supported by I.N.d.A.M. while preparing this work.
He would like to thank F. Bonechi, R. Giachetti, E. Sorace
and M. Tarlini for fruitful conversations on the subject.

\enddocument